\documentclass[fleqn,11pt]{article}

\usepackage{amsmath,amsfonts,amssymb}
\usepackage{color}
\usepackage{geometry}
\newtheorem{theorem}{Theorem}
\newtheorem{lemma}[theorem]{Lemma}

\newtheorem{proposition}{Proposition}

\newcommand{\R}{\mathbb{R}}

\newcommand{\N}{\mathbb{N}}

\newcommand{\cqfd}
{%
\mbox{}%
\nolinebreak%
\hfill%
\rule{2mm}{2mm}%
\newline
\newline
}

\title{
Discrete velocity Boltzmann equations in the plane: stationary solutions.}
\author{Leif ARKERYD and Anne NOURI\\
\\Mathematical Sciences, 41296 G\"oteborg, Sweden,\\
arkeryd@chalmers.se\\
\hspace*{1.in}\\
Aix-Marseille University, CNRS, Centrale Marseille, I2M UMR 7373, 13453 Marseille, France,\\
anne.nouri@univ-amu.fr\\
www.i2m.univ-amu.fr/perso/anne.nouri/}

\date{}

\begin{document}
\maketitle
\hspace{1cm}\\
{\bf Abstract} \hspace{1cm}\\
The paper proves existence of stationary mild solutions for normal discrete velocity Boltzmann equations in the plane with no pair of colinear interacting velocities and given ingoing boundary values. An important restriction of all velocities pointing into the same half-space in a previous paper is removed in this paper. A key property is $L^1$ compactness of integrated collision frequency for a sequence of approximations. This is proven using the Kolmogorov-Riesz theorem, which here replaces the $L^1$ compactness of velocity averages in the continuous velocity case, not available when  the velocities are discrete.
\footnotetext[1] {2010 Mathematics Subject Classification; 60K35, 82C40, 82C99.} \\
\footnotetext[2] {Key words; stationary Boltzmann equation, discrete coplanar velocities, normal model.}.
%
% SECTION 1
%
\section {Introduction.}
\label{generalization}
\setcounter{equation}{0}
The Boltzmann equation is the fundamental mathematical model 
in the kinetic theory of gases. Replacing its continuum of velocities with a discrete set of velocities is a  simplification,  preserving the  essential features of free flow and quadratic collision term. Besides this fundamental aspect, the discrete equations can  approximate the Boltzmann equation with any given accuracy \cite {BPS}, \cite{FKW}, \cite{M}, and are thereby useful for approximations and numerics. In the quantum realm they can also be more directly
connected to microscopic quasi/particle models.
A discrete velocity model of a kinetic gas is a system of partial differential equations having the form,
\begin{align*}
&\frac{\partial f_i}{\partial t}(t,z)+v_i\cdot \nabla _zf_i(t,z)= Q_i(f,f)(t,z),\quad t>0,\quad z\in \Omega ,\quad 1\leq i\leq p,
\end{align*}
where $f_i(t,z)$, $1\leq i\leq p$, are phase space densities at time $t$, position $z$ and velocities $v_i$. The spatial domain is $\Omega %\subset \R ^2
$. The given discrete velocities are $v_i$, $1\leq i\leq p$. For $f= (f_i)_{1\leq i\leq p}$, the collision operator $Q= (Q_i)_{1\leq i\leq p}$
with gain part $Q^+$, loss part $Q^-$, and collision frequency $\nu$, is given by
\begin{align*}
&Q_i(f,f)= \sum_{j,l,m=1}^p \Gamma_{ij}^{lm} (f_lf_m-f_if_j)\\
&\hspace*{0.52in}= Q_i^+(f,f)-Q_i^-(f,f), \\
&Q^+_i(f,f)= \sum_{j,l,m=1}^p \Gamma_{ij}^{lm}f_lf_m, \quad Q^-_i(f,f)= f_i\nu_i(f), \quad \nu _i(f)= \sum_{j,l,m=1}^p \Gamma_{ij}^{lm} f_j,\quad i=1,..., p.
\end{align*}
The collision coefficients satisfy
\begin{align}\label{Gamma}
&\Gamma_{ij}^{lm}=\Gamma_{ji}^{lm}=\Gamma_{lm}^{ij}\geq 0.
\end{align}
If a  collision coefficient $\Gamma_{ij}^{lm}$ is non-zero, then the conservation laws for momentum and energy,
\begin{align}\label{conservations-velocities}
&v_i+v_j=v_l+v_m,\quad |v_i|^2+|v_j|^2=|v_l|^2+|v_m|^2,
\end{align}
are satisfied. We call interacting velocities any couple of velocities $(v_i,v_j)$ such that for some $(l,m)\in \{ 1,\cdot \cdot \cdot ,p\} ^2 $, $\Gamma _{ij}^{lm}> 0$. The discrete velocity model (DVM) is called normal (see \cite{C}) if any solution of the equations
\begin{eqnarray*}
\Psi(v_i)+\Psi(v_j)=\Psi(v_l)+\Psi(v_m),
\end{eqnarray*}
where the indices $(i,j;l,m)$ take all possible values satisfying $\Gamma_{ij}^{lm}>0$, is given by
\begin{align*}
&\Psi(v)=a+b\cdot v+c|v|^2,
\end{align*}
for some constants $a,c\in \R$ and $b\in \R ^d$. 
We consider
\begin{align}
&\text{the generic case of normal coplanar velocity sets with}\nonumber \\
&\text{ no pair of  colinear interacting velocities  } (v_i, v_j).\label{generic-condition}
\end{align}
The case is generic. Indeed, consider a normal velocity set such that for some interacting velocities $(v_i,v_j)$, $v_i$ and $v_j$ are colinear. Then there exists an arbitrary small vector $v_0$ such that the velocity set $(v_i+v_0)_{1\leq i\leq p} $ is normal and with no colinear interacting velocities.
The paper considers stationary solutions to normal coplanar discrete velocity models satisfying \eqref{generic-condition}, in a strictly convex bounded open subset $\Omega \subset \R ^2$, with  $C^2$ boundary $\partial \Omega $ and given boundary inflow. 
Denote by $n(Z)$ the inward normal to $Z\in \partial \Omega $. Denote the $v_i$-ingoing (resp. $v_i$-outgoing) part of the boundary by
\begin{align*}
&\partial \Omega _i^+= \{ Z\in \partial \Omega \hspace*{0.02in}; \hspace*{0.02in}v_i\cdot n(Z)>0\}  ,\quad \quad (\text{resp.   }\partial \Omega _i^-= \{ Z\in \partial \Omega \hspace*{0.02in}; \hspace*{0.02in}v_i\cdot n(Z)<0\}  ).
\end{align*}
Let
\begin{align*}
&s_i^+(z)= \inf \{ s>0\hspace*{0.02in}; z-sv_i\in \partial \Omega _i^+\} ,\quad  s_i^-(z)= \inf \{ s>0\hspace*{0.02in}; z+sv_i\in \partial \Omega _i^-\} ,\quad  z\in \Omega .
\end{align*}
Write
\begin{align}\label{df-in-out-points}
z_i^+(z)= z-s_i^+(z)v_i \quad (\text{resp.   } z_i^-(z)= z+s_i^-(z)v_i)
\end{align}
for the ingoing (resp. outgoing) point on $\partial\Omega$ of the characteristics through $z$ in direction $v_i$.\\
The stationary boundary value problem
\begin{align}
&v_i\cdot \nabla f_i(z)= Q_i(f,f)(z),\quad z\in \Omega,\label{discrete-general-a}\\
& f_i(z)= f_{bi}(z),\quad z\in \partial \Omega _i^+,\quad \quad 1\leq i\leq p,\label{discrete-general-b}
\end{align}
is considered in $L^1$ in one of the following equivalent
%(cf \cite{DPL})(where the proof of equivalence is practically the same as their last five lines of page 362 and page 363 )
forms (\cite{DPL});\\
the exponential multiplier form,
\begin{align}\label{exponential-form}
f_{i}(z)&= f_{bi}(z_i^+(z))e^{-\int_{0}^{s_i^+(z)} \nu_i(f)(z_i^+(z)+sv_i)ds}\nonumber\\
&+\int_{0}^{s_i^+(z)} Q_i^+(f,f)(z_i^+(z)+sv_i)e^{-\int_s^{s_i^+(z)}\nu_i(f)(z_i^+(z)+rv_i)dr}ds ,\quad \text{a.a.  }z\in \Omega ,\quad 1\leq i\leq p,
\end{align}
the mild form,
\begin{align}\label{mild-form}
&f_i(z)= f_{bi}(z_i^+(z))+\int_{0}^{s_i^+(z)} Q_i(f,f)(z_i^+(z)+sv_i)ds,\quad \text{a.a.  }z\in \Omega ,\quad 1\leq i\leq p,
\end{align}
the renormalized form,
\begin{align}\label{renormalized-form}
v_i\cdot \nabla \ln(1+f_i)(z)= \frac{Q_i(f,f)}{1+f_i}(z),\hspace*{0.04in} z\in \Omega ,\quad \quad \quad f_i(z )= f_{bi}(z),\hspace*{0.04in} z\in \partial \Omega _i^+,\quad 1\leq i\leq p,
\end{align}
in  the sense of distributions.
Denote by $L^1_+(\Omega)$ the set of non-negative integrable functions on $\Omega$. For a distribution function $f= (f_i)_{1\leq i\leq p}$, define its entropy (resp. entropy dissipation)  by
\begin{align*}
&\sum _{i= 1}^p\int _\Omega f_i\ln f_i(z)dz,\quad \Big( \text{resp.}\quad \sum_{i,j,l,m=1}^p \Gamma_{ij}^{lm} \int _\Omega (f_lf_m-f_if_j)\ln \frac{f_lf_m}{f_if_j}(z)dz\Big) .
\end{align*}
The main result of the paper is
%
% Theorem 1.1
%
\begin{theorem}\label{main-result}
\hspace*{1.in}\\
Consider a coplanar normal discrete velocity model and a non-negative ingoing boundary value $f_{b}$ with mass and entropy inflows bounded,
\begin{align*}
&\int _{\partial \Omega _i^+}v_i\cdot n(z)\hspace*{0.02in}f_{bi}(1+\ln f_{bi})(z)d\sigma (z)<+\infty,\quad 1\leq i\leq p.
\end{align*}
 For the boundary value problem \eqref{discrete-general-a}-\eqref{discrete-general-b} satisfying \eqref{generic-condition}, there exists a stationary mild solution in $\big( L^1_+(\Omega )\big) ^p$  with finite mass and entropy-dissipation.
\end{theorem}
\hspace*{1.in}\\
Given $i\in \{ 1,\cdot \cdot \cdot , p\} $, if $\Gamma _{ij}^{lm}= 0$ for all $j$, $l$ and $m$, then $f_i$ equals its ingoing boundary value, and the rest of the system can be solved separately. Such $i$'s are not present in the following discussion.\\
Most mathematical results for stationary discrete velocity models of the Boltzmann equation have been obtained in one space dimension. An overview is given in \cite{IP}. Half-space problems \cite{NB1} and weak shock waves \cite{NB2} for discrete velocity models have also been studied. A discussion of normal discrete velocity models, i.e. conserving nothing but mass, momentum and energy, can be found in \cite{BVW}.
In two dimensions, special classes of solutions to the Broadwell model are given in \cite{BT}, \cite{B} and  \cite{Ily}. The Broadwell model, not included in the present results, is a four-velocity model, with $v_1+v_2=v_3+v_4=0$ and $v_1$, $v_3$ orthogonal. \cite{B} contains a detailed study of the stationary Broadwell equation in a rectangle with comparison to a
Carleman-like system, and a discussion of (in)compressibility aspects. A main result in \cite{CIS} is the existence of continuous solutions to the two-dimensional stationary Broadwell model with continuous boundary data for a rectangle. The paper \cite{AN2} solves that problem in an $L^1$-setting. The proof uses in an essential way the constancy of the sums  $f_1+f_2$ and $f_3+f_4$ along characteristics, which no longer holds in the present paper. For every normal model, there is a priori control of entropy dissipation, mass and entropy flows through the boundary. 
From there, main difficulties are to prove that for a sequence of approximations, weak $L^1$ compactness  holds and the limit of the collision operator equals the collision operator of the limit. In \cite{AN1}, weak $L^1$ compactness of a sequence of approximations was obtained with assumption \eqref{generic-condition} together with the assumption that all velocities $v_i$ point out into the same half-plane. In this paper we keep assumption \eqref{generic-condition}, remove the second assumption and provide a new proof of weak $L^1$ compactness of approximations using \eqref{generic-condition}. Assumption \eqref{generic-condition} is also crucial for proving $L^1$ compactness of the integrated collision frequencies, that is important for the convergence procedure. Our paper also differs from \cite{AN1} in the limit procedure. The frame of the limit procedure in \cite{AN1} is the splitting into 'good' and 'bad' characteristics following the approach in our earlier stationary continuous velocity papers \cite{AN}-\cite{AN4}. Here we have instead recourse to sub- and super-solutions used in the classical evolutionary frame for renormalized solutions to the Boltzmann equation \cite{DPL}.
\\
\hspace*{1.in}\\
For the continuous velocity evolutionary Boltzmann equation \cite{DPL}, the compactness properties of the collision frequency use in an essential way the averaging lemma, which is not available for the discrete velocity Boltzmann model. In the present paper, the compactness properties are proven by the Kolmogorov-Riesz theorem.
Also the argument used in the stationary paper \cite{AN4} in the continuous velocity case for obtaining control of entropy, hence weak $L^1$ compactness of a sequence of approximations  from the control of entropy dissipation, does not work in a discrete velocity case because the number of velocities is finite.\\
\hspace*{1.in}\\
The proof starts in Section \ref{Section2} from bounded  approximations. In Section \ref{Section3}, $L^1$ compactness properties of the approximations are proven. Section \ref{Section4} is devoted to the proof of Theorem \ref{main-result}.
%
%
% SECTION 2
%
\section{Approximations.}\label{Section2}
\label{approximations}
\setcounter{equation}{0}
\setcounter{theorem}{0}
Denote by $\N ^*= \N \setminus \{ 0\} $ and by $a\wedge b$ the minimum  of two real numbers $a$ and $b$. Let $\mu_\alpha$ be a smooth mollifier in $\R^2$ with support in the ball centered at the origin of radius $\alpha $. Outside the boundary the function to be convolved with 
$\mu_\alpha$\textcolor{magenta}{,}
is continued in the normal direction by its boundary value. Let $\tilde{\mu }_{k}$ be a smooth mollifier on $\partial \Omega $ in a ball of radius $\frac{1}{k}$. Denote by
\begin{align*}
&f^k_{bi}= \Big( f_{bi}(\cdot )\wedge \frac{k}{2}\Big) \ast \tilde{\mu }_k ,\quad 1\leq i\leq p, \quad k\in \N ^*.
\end{align*}
The lemma introduces a primary approximated boundary value problem with damping and convolutions.
\begin{lemma}\label{first-approximations}
\hspace*{1.in}\\
For any $\alpha >0$ and $k\in \N ^*$, there is a solution $F^{\alpha ,k}\in (L^1_+(\Omega))^p$ to
\begin{align}
&\alpha F^{\alpha ,k}_i+v_i\cdot \nabla F^{\alpha ,k}_i= \sum _{j,l,m=1}^p\Gamma _{ij}^{lm}\Big( \frac{F^{\alpha ,k}_l}{1+\frac{F^{\alpha ,k}_l}{k}}\frac{F^{\alpha ,k}_m\ast \mu_\alpha }{1+\frac{F^{\alpha ,k}_m\ast \mu_\alpha }{k}}-\frac{F^{\alpha ,k}_i}{1+\frac{F^{\alpha ,k}_i}{k}}\frac{F^{\alpha ,k}_j\ast \mu_\alpha }{1+\frac{F^{\alpha ,k}_j\ast \mu_\alpha }{k}}\Big)\hspace*{0.03in},\label{df-F-final-1}\\
&F^{\alpha ,k}_i(z)= f^k_{bi}(z),\quad z\in \partial \Omega _i^+,\quad \quad 1\leq i\leq p.\label{df-F-final-2}
\end{align}
\end{lemma}
\underline{Proof of Lemma \ref{first-approximations}.}\\
For a proof of Lemma \ref{first-approximations} we refer to the second section in \cite{AN1}.\\
\hspace*{1.in}\\
Let $k\in \N ^*$ be given. Each component of $F^{\alpha ,k}$ is bounded by a multiple of $k^2$. Therefore $(F^{\alpha ,k})_{\alpha \in ]0,1[}$ is weakly compact in $(L^1(\Omega))^p$.  For a subsequence, the convergence is strong in $(L^1(\Omega))^p$ as stated in the following lemma.
%
% Lemma 2,2
%
\begin{lemma}\label{cv-F-alpha}
\hspace*{1.in}\\
{There is a sequence $(\beta (q))_{q\in \N}$ tending to zero when $q\rightarrow +\infty $ and a function $F^k\in L^1$,  such that $(F^{\beta (q),k})_{q\in \N }$ } strongly converges in $(L^1(\Omega))^p$ to $F^k$ when $q\rightarrow +\infty $.
\end{lemma}
\underline{Proof of Lemma \ref{cv-F-alpha}.}\\
For a proof of Lemma \ref{cv-F-alpha} we refer to Lemma 3.1 in \cite{AN1}.\\
\hspace*{1.in}\\
Denote by
\begin{align}\label{df-approximate-gain}
&Q^{+k}_i= \sum _{j,l,m=1}^p\Gamma _{ij}^{lm}\frac{F^k_l}{1+\frac{F^k_l}{k}}\frac{F^k_m}{1+\frac{F^k_m}{k}},\quad \nu _i^k= \sum _{j,l,m=1}^p\Gamma _{ij}^{lm}\frac{F^k_j}{(1+\frac{F_i^k}{k})(1+\frac{F_j^k}{k})},\nonumber \\
&Q_i^k= Q_i^{+k}-F_i^k\nu _i^k,\quad 1\leq i\leq p,
\end{align}
and by $\tilde{D}_k$ the entropy production term of the approximations,
\begin{align}\label{df-entropy-production}
&\tilde{D}_k= \sum _{i, j,l,m=1}^p\Gamma _{ij}^{lm}\Big( \frac{F^k_l}{1+\frac{F^k_l}{k}}\frac{F^k_m}{1+\frac{F^k_m}{k}}-\frac{F^k_i}{1+\frac{F^k_i}{k}}\frac{F^k_j}{1+\frac{F^k_j}{k}}\Big) \ln \frac{F^k_lF^k_m(1+\frac{F^k_i}{k})(1+\frac{F^k_j}{k})}
{(1+\frac{F^k_l}{k})(1+\frac{F^k_m}{k})F^k_iF^k_j}\hspace*{0.04in}.
\end{align}
All along the paper, $c_b$ denotes constants that may vary from line to line but is independent of parameters tending to $+\infty $ or to zero.
%
%
% Lemma 2.3
%
\begin{lemma}\label{existence-k-approximations}
\hspace*{1.in}\\
$F^k$ is a non-negative solution to 
\begin{align}
&v_i\cdot \nabla F_i^k=Q_i^{+k}-F_i^k\nu_i^k
\hspace*{0.03in},\label{Fk-i}\\
&F^k_i(z)= f^k_{bi}(z),\quad z\in \partial \Omega _i^+,\quad 1\leq i\leq p.\label{bcFk-i}
\end{align}
Solutions $(F^k)_{k\in \N^*}$ to \eqref{Fk-i}-\eqref{bcFk-i} have mass and entropy dissipation bounded from above uniformly with respect to $k$. 
Moreover their outgoing flows at the boundary are controlled as follows,
\begin{align}\label{outgoing-flows}
&\sum_{i= 1}^p\int_{\partial \Omega _i^{-},F_i^k\leq k}  \mid v_i\cdot n(Z)\mid F^k_i\ln F^k_i(Z)d\sigma (Z)+\ln \frac{k}{2}\int_{\partial \Omega _i^{-},F_i^k> k}\mid v_i\cdot n(Z)\mid {F^k_i}d\sigma (Z)\leq c_b.
\end{align}
\end{lemma}
\underline{Proof of Lemma \ref{existence-k-approximations}.}\\%\ref{existence-k-approximations}}\\
Passing to the limit when $q\rightarrow +\infty $ in \eqref{df-F-final-1}-\eqref{df-F-final-2} written for $F^{\beta (q),k}$, implies that $F^k$ is a solution in $\big( L^1_+(\Omega )\big) ^p$ to \eqref{Fk-i}-\eqref{bcFk-i}. For a proof of the rest of Lemma \ref{existence-k-approximations}, we refer to Lemma 3.2 in \cite{AN1}.\\
%
%
%
%
%
% SECTION 3
%
\section{On compactness of sequences of approximations.}\label{Section3}
\setcounter{equation}{0}
\setcounter{theorem}{0}

This section is devoted to prove $L^1$ compactness properties of the approximations. In Proposition \ref{weakL1compactness}, weak $L^1$ compactness of $(F^k)_{k\in \N ^*}$ is proven. Lemma \ref{lemma-df-Omega-epsilon} splits $\Omega $ into a set of $i$-characteristics with arbitrary small measure and its 
complement\textcolor{magenta}{,} where both the approximations and their integrated collision frequencies are bounded. In Lemma \ref{compactness-integrated-collision-frequency}, 
the strong $L^1$ compactness of integrated collision frequency is proven.
%
% Proposition 3.1
%
 \begin{proposition}\label{weakL1compactness}
 \hspace*{1.in}\\
 The sequence $(F^k)_{k\in \N ^*}$ solution to \eqref{Fk-i}-\eqref{bcFk-i} is weakly compact in $L^1$. 
 \end{proposition}
\underline{Proof of Proposition \ref{weakL1compactness}.}\hspace*{0.05in}\\
By Lemma \ref{existence-k-approximations}, $(F^k)_{k\in \N ^*}$ is uniformly bounded in $(L^1(\Omega)) ^p$. \\
\hspace*{1.in}\\
Given \eqref{outgoing-flows} and the following bound on $F^k$,
\begin{align}\label{bdd-by-outgoingA}
F^k_i(z)&\leq F^k_i(z+s_i^-(z)v_i)\hspace*{0.02in}\exp \Big( \Gamma \displaystyle \sum _{j\in J_i}\int _{-s_i^+(z)}^{s_i^-(z)}F_j(z+rv_i)dr\Big) ,\quad z\in \Omega ,\quad i\in \{ 1,\cdot \cdot \cdot ,p\} ,
\end{align}
the weak $L^1$ compactness of $(F^k)_{k\in \N ^*}$ will follow from the uniform boundedness in $L^\infty (\partial \Omega _i^+)$ of 
\begin{align}\label{weak29novB}
&\big( \int _{0}^{s_i^-(Z)}F_j(Z+rv_i)dr\big) _{j\in J_i, k\in \N },
\end{align}
where $J_i$ denotes $\{ j\in \{ 1,\cdot \cdot \cdot , p\} ; (v_i,v_j)\text{   are interacting velocities}\} $. By \eqref{generic-condition}, there exists $\eta >0$ such that for all interacting velocities $(v_i,v_j)$,
\begin{align}\label{orthogonal-assumption}
&\lvert \sin (\widehat{v_i,v_j})\rvert >\eta .
\end{align}
Let $i\in \{ 1,\cdot \cdot \cdot ,p\} $ and $Z\in \partial \Omega _i^+$. Multiply the equation satisfied by $F_j^k$ by $\frac{v_i^\perp \cdot v_j}{\lvert v_i\rvert }$ and integrate it on one of the half domains defined by the segment $[Z,Z+s_i^-(Z)v_i] $. Summing over $j\in \{ 1,\cdot \cdot \cdot ,p\}$ implies that
\begin{align}\label{weak29novA}
&\sum _{j=1}^p\sin ^2(\widehat{v_i,v_j})\int _0^{s_i^-(Z)}F_j^k(Z+sv_i)ds\leq c_b,\quad Z\in \partial \Omega _i^+.
\end{align}
Together with \eqref{orthogonal-assumption}, this leads to the control of \eqref{weak29novB}.  \cqfd 
\hspace*{1.in}\\
Recall the exponential multiplier form for the approximations $(F^k)_{k\in \N ^*}$,
\begin{align}\label{exponential-form-approximations}
F^k_{i}(z)&= f^k_{bi}(z_i^+(z))e^{-\int_{-s_i^+(z)}^0 \nu_i^k(z+sv_i)ds}\nonumber\\
&+\int_{-s_i^+(z)}^0 Q_i^{+k}(z+sv_i)e^{-\int_s^{0}\nu_i^k(F^k)(z+rv_i)dr}ds ,\quad \text{a.a.  }z\in \Omega ,\quad 1\leq i\leq p,
\end{align}
with $\nu _i^k$ and $Q_i^{+k}$ defined in \eqref{df-approximate-gain}. An $i$-characteristics is a segment of points $[Z-s_i^+(Z)v_i,Z]$, where $Z\in \partial\Omega _i^-$. 
Denote by $\Gamma = \displaystyle \max _{i,j,l,m}\Gamma _{ij}^{lm}$.
%
%  Lemma 3.1
%
\begin{lemma}\label{lemma-df-Omega-epsilon}
\hspace*{1.in}\\
For $i\in \{ 1, ..., p\} $, $k\in \N ^*$ and $\epsilon >0$, there is a subset $\Omega ^{k,\epsilon }_{i}$ of $i$-characteristics of $\Omega $ with measure smaller than $c_b\epsilon $,
such that for any $z\in \Omega \setminus \Omega ^{k,\epsilon }_{i}$,
\begin{align}\label{bounds-on-complementary-of-Omega_iEpsilon}
&F^k_{i}(z)\leq \frac{1}{\epsilon ^2}\exp \big( \frac{p\Gamma }{\epsilon ^2}\big) ,\quad \int_{-s_i^+(z)}^{s_i^-(z)}\nu_i^k(z+sv_i)ds\leq \frac{p\Gamma }{\epsilon ^2}.
\end{align}
\end{lemma}
\underline{Proof of Lemma \ref{lemma-df-Omega-epsilon}.}\\
By the strict convexity of $\Omega $, there are for every $i\in \{ 1,\cdot \cdot \cdot p\} $ two points of $\partial \Omega $, denoted by $\tilde{Z_i}$ and $\bar{Z_i}$ such that
\begin{align*}
v_i\cdot n(\tilde{Z_i})= v_i\cdot n(\bar{Z_i})= 0.
\end{align*}
Let $\tilde{l}_i$ (resp. $\bar{l}_i$) be the largest boundary arc included in $\partial \Omega _i^-$ with one end point $\tilde{Z}_i$ (resp. $\bar{Z}_i$) such that
\begin{align}\label{lemma4-1b}
-\epsilon \leq v_i\cdot n(Z)\leq 0,\quad Z\in \tilde{l}_i\cup \bar{l}_i.
\end{align}
Let $J_i$ be the subset of $\{ 1,\cdot \cdot \cdot , p\} $ such that 
\begin{align}\label{df-J-i}
&\text{for some  }(l,m)\in \{ 1, \cdot \cdot \cdot ,p \} ^2,\quad \Gamma _{ij}^{lm}>0, \quad j\in J_i.
\end{align} 
It follows from the exponential form of $F^k_i$ that
\begin{align}\label{bdd-by-outgoing}
F^k_i(z)&\leq F^k_i(z+s_i^-(z)v_i)\hspace*{0.02in}\exp \Big( \Gamma \displaystyle \sum _{j\in J_i}\int _{-s_i^+(z)}^{s_i^-(z)}F_j(z+rv_i)dr\Big) ,\quad z\in \Omega .
\end{align}
The boundedness of the mass flow of $(F^k_i)_{k\in \N ^*}$ across $\partial \Omega _i^-$ is
\begin{align}\label{lemma4-1a}
&\int _{\partial \Omega _i^-}\mid v_i\cdot n(Z)\mid F^k_i(Z)d\sigma (Z)\leq c_b,\quad k\in \N ^*.
\end{align}
It follows from \eqref{lemma4-1b}-\eqref{lemma4-1a} that the measure of the set
\begin{align*}
\{ Z\in \partial \Omega _i^-\cap \tilde{l}_i^{\hspace*{0.02in}c}\cap \bar{l}_i^{\hspace*{0.02in}c}\quad ;\quad  F^k_i(Z)>\frac{1}{\epsilon ^2}\}
\end{align*}
is smaller than $c_b\epsilon $. The boundedness of the mass of $(F^k_j)_{k\in \N ^*}$ can be written
\begin{align*}
&\int_\Omega F_j^k(z)dz=\int _{\partial \Omega _i^-}\mid v_i\cdot n(Z)\mid \Big( \int _{-s_i^+(Z)}^0F^k_j(Z+rv_i)dr\Big) d\sigma (Z)\leq c_b,\quad j\in J_i.
\end{align*}
Hence the measure of the set
\begin{align*}
&\{ Z\in \partial \Omega _i^-\cap \tilde{l}_i^{\hspace*{0.02in}c}\cap \bar{l}_i^{\hspace*{0.02in}c}\quad ;\quad  \int _{-s_i^+
(Z)}^0F^k_j(Z+rv_i)dr>\frac{1}{\epsilon ^2}\} ,\quad j\in J_i,
\end{align*}
is smaller than $c_b\epsilon $. Consequently, the measure of the set of $Z\in \partial \Omega _i^-\cap \tilde{l}_i^{\hspace*{0.02in}c}\cap \bar{l}_i^{\hspace*{0.02in}c}$ outside of which
\begin{align*}
F^k_i(Z)\leq\frac{1}{\epsilon ^2} \quad  \text{and} \quad \int _{-s_i^+(Z)}^0F^k_j(Z+rv_i)dr\leq\frac{1}{\epsilon ^2},\quad j\in J_i,
\end{align*}
is bounded by $c_b\epsilon$.
Together with \eqref{bdd-by-outgoing}, this implies that the measure of the complement of the set of
$Z\in \partial \Omega _i^-$, such that
\begin{align*}
&F_i^k(z)\leq \frac{1}{\epsilon ^2}\exp \big( \frac{p\Gamma }{\epsilon ^2}\big) \quad \text{and}\quad \int _{-s_i^+(z)}^{s_i^-(z)}\nu _i^k(z+rv_i)dr\leq \frac{p\Gamma }{\epsilon ^2}
\end{align*}
for $z=Z+sv_i$, $s\in [ -s_i^+(Z),0]$, is bounded by $c_b\epsilon$. With it $c_b\epsilon$ is a bound for the measure of the complement, denoted by $\Omega _{i}^{k,\epsilon }$, of the set of $i$-characteristics in $\Omega$ such that for all points $z$ on the $i$-characteristics, \eqref{bounds-on-complementary-of-Omega_iEpsilon} holds.
\cqfd
Given $i\in \{ 1, ..., p\} $ and $\epsilon >0$, let $\chi^{k,\epsilon }_{i}$ denote the characteristic  function of the complement of $\Omega ^{k,\epsilon }_{i}$. The following lemma proves the compactness in $L^1(\Omega )$ of the $k$-sequence of integrated collision frequencies.
%
% Lemma 3.2
%
\begin{lemma}\label{compactness-integrated-collision-frequency}
\hspace*{1.in}\\
The sequences $\Big( \int _{-s_i^+(z)}^0 \nu _i^k(z+sv_i)ds\Big) _{k\in \N ^*}$, $1\leq i\leq p$, are strongly compact in $L^1(\Omega )$.
\end{lemma}
\underline{Proof of Lemma \ref{compactness-integrated-collision-frequency}.}\\
Take $\Gamma_{ij}^{lm}> 0$. By \eqref{generic-condition}, $v_i$ and
$v_j$ span $\R^2$. Denote by $(a,b)$ the corresponding coordinate system, $(a^-, a^+)$ defined by
\begin{align*}
&a^-= \min \{ a\in \R ; (a,b)\in \Omega \hspace*{0.04in}\text{for some   }b\} ,\quad a^+= \max \{ a\in \R ; (a,b)\in \Omega \hspace*{0.04in}\text{for some   }b\} ,
\end{align*}
and by $D$ the Jacobian of the change of variables $z\rightarrow (a,b)$.
The uniform bound for the mass of $(F^k)_{k\in \N ^*}$ proven in Lemma \ref{existence-k-approximations}, implies that
\begin{align*}
&\Big( \int _\Omega \int _{-s_i^+(z)}^0 \nu _i^k(z+sv_i)dsdz\Big) _{k\in \N ^*}
\end{align*}
is bounded in $L^1$ uniformly with respect to $k$. Indeed, for some $(b^-(a),b^+(a))$, $a\in [a^-,a^+]$,
\begin{align*}
\int _\Omega \int _{-s_i^+(z)}^0 F_j^k(z+sv_i)dsdz&= D\int _{a^-}^{a^+}\int _{b^-(a)}^{b^+(a)}\int _{-s_i^+(bv_j)}^{a}F_j^k(bv_j+sv_i) ds\hspace*{0.02in}db\hspace*{0.02in}da\\
&\leq D\int _{a^-}^{a^+}\int _{b^-(a)}^{b^+(a)}\int _{-s_i^+(bv_j)}^{s_i^-(bv_j)}F_j^k(bv_j+sv_i) ds\hspace*{0.02in}db\hspace*{0.02in}da\\
&\leq c\int _\Omega F_j^k(z)dz,\hspace*{1.4in} j\in J_i.
\end{align*}
By the Kolmogorov-Riesz theorem (\cite{K}, \cite{R}), the compactness of $\Big( \int _{-s_i^+(z)}^0 \nu _i^k(z+sv_i)ds\Big) _{k\in \N ^*}$ will follow from its translational equi-continuity in $L^1(\Omega)$. Equicontinuity in the direction $v_i$, and in the direction $v_j$ with the mild form \eqref{mild-form} for $F_j^k$, come natural. Here the assumption \eqref{generic-condition} becomes crucial.
The sequence
\begin{align}\label{pf-lemma4-2-a}
&\Big( \int _{-s_i^+(z)}^0 F_j^k(z+sv_i)ds\Big) _{k\in \N ^*},\quad j\in J_i,
\end{align}
is translationally equi-continuous in the $v_i$-direction. Indeed, $s_i^+(z+hv_i)= s_i^+(z)+h$
so that, denoting by $I(0,h)$ the interval with endpoints $0$ and $h$ and using the uniform bound on the mass of $(F_j^k)_{k\in \N ^*}$,
\begin{align*}
&\int _\Omega \mid \int _{-s_i^+(z+hv_i)}^0 F_j^k(z+hv_i+sv_i)ds-\int _{-s_i^+(z)}^0 F_j^k(z+sv_i)ds\mid dz\\
&= \int _\Omega \int _{s\hspace*{0.02in}\in I(0,h)}F_j^k(z+sv_i)ds\hspace*{0.02in}dz\\
&\leq c \mid h\mid .
\end{align*}
Let us prove the translational equi-continuity of \eqref{pf-lemma4-2-a} in the $v_j$-direction. By the weak $L^1$ compactness of $(F_j^k) _{k\in \N ^*}$, it is sufficient to prove the translational equi-continuity in the $v_j$-direction of $\big( \int _{s_i^+(z)}^0 \chi _j^{k,\epsilon }F_j^k(z+sv_i)ds\big) _{k\in \N ^*}$.
Expressing $F_j^k(z+hv_j+sv_i)$ (resp. $F_j^k(z+sv_i)$) as integral along its $v_j$-characteristics, it holds that
\begin{align*}
&\mid \int _{-s_i^+(z+hv_j)}^0 \chi _j^{k,\epsilon }F_j^k(z+hv_j+sv_i)ds-\int _{-s_i^+(z)}^0 \chi _j^{k,\epsilon }F_j^k(z+sv_i)ds\mid \leq \mid A_{ij}^k(z,h)\mid +\mid B_{ij}^k(z,h)\mid ,
\end{align*}
where
\begin{align*}
&A_{ij}^k(z,h)= \int _{-s_i^+(z+hv_j)}^0 \chi _j^{k,\epsilon }f_{bj}^k\big( z_j^+(z+hv_j+sv_i) \big) ds-\int _{-s_i^+(z)}^0 \chi _j^{k,\epsilon }f_{bj}^k\big( z_j^+(z+sv_i)\big) ds,
\end{align*}
and
\begin{align*}
B_{ij}^k(z,h)&= \int _{-s_i^+(z+hv_j)}^0 \int _{-s_j^+(z+hv_j+sv_i)}^0\chi _j^{k,\epsilon }Q_j^k(z+hv_j+sv_i+rv_j)drds\nonumber \\
&-\int _{-s_i^+(z)}^0 \int _{-s_j^+(z+sv_i)}^0\chi _j^{k,\epsilon }Q_j^k(z+sv_i+rv_j)drds,
\end{align*}
with $Q_i^k$ defined in \eqref{df-approximate-gain}. Denote by $(z_j^+(z_i^+(z)), z_j^+(z_i^+(z+hv_j))$ the boundary arc with end points $z_j^+(z_i^+(z))$ and 
$z_j^+(z_i^+(z+hv_j))$ and of length tending to zero with $h$. Performing the change of variables $s\rightarrow Z= z_j^+(z+hv_j+sv_i) $ (resp. $s\rightarrow Z= z_j^+(z+sv_i) $) in the first (resp. second) term of $A_{ij}^k(z,h)$, and using that the sequence $(f_{bi}^k)_{k\in \N ^*}$ is bounded by $f_{bi}$, it holds that 
\begin{align}\label{pf-lemma4-2-e}
&\lim _{h\rightarrow 0}\int _\Omega \mid A_{ij}^k(z,h)\mid dz= 0,
\end{align}
uniformly with respect to $k$. Moreover, for some $\omega _h(z)\subset \Omega $ of measure or order $\mid h\mid $ uniformly with respect to $z\in \Omega $,
\begin{align}\label{pf-lemma4-2-b}
B_{ij}^k(z,h)&= \int _{\omega _h(z)} \chi _j^{k,\epsilon }Q_j^k(Z)dZ.
\end{align}
The sequence $(\chi _j^{k,\epsilon }Q_j^k)_{k\in \N ^*}$ is weakly compact in $L^1$. Indeed,
\begin{align}\label{pf-lemma4-2-c}
\chi _j^{k,\epsilon }Q_j^k&\leq \frac{1}{\ln \Lambda }\tilde{D}_k+\Gamma \Lambda \Big( \sum _{i\in J_j}F_i^k\Big) (\chi _j^{k,\epsilon }F_j^k)\nonumber \\
&\leq \frac{1}{\ln \Lambda }\tilde{D}_k+\frac{\Gamma \Lambda }{\epsilon ^2}\exp \big( \frac{p\Gamma }{\epsilon ^2}\big) \Big( \sum _{i\in J_j}F_i^k\Big) ,\quad \Lambda >1,
\end{align}
with $(\tilde{D}_k)_{k\in \N ^*}$ uniformly bounded in $L^1$ and $(F_i^k)_{k\in \N ^*}$ weakly compact in $L^1$. Hence,
\begin{align}\label{pf-lemma4-2-g}
&\lim _{h\rightarrow 0}\int _\Omega \mid B_{ij}^k(z,h)\mid dz= 0,\quad \text{uniformly with respect to  }k.
\end{align}
\cqfd
%
%
%
%
% SECTION 4
%
%
%
%
\section{The passage to the limit in the approximations.}\label{Section4}
\setcounter{equation}{0}
\setcounter{theorem}{0}
%
%
% Definition of the final test function chi _eta
%
Let $f$ be the weak $L^1$ limit of a subsequence of the solutions $(F^k)_{k\in \N ^*}$ to \eqref{Fk-i}-\eqref{bcFk-i}, still denoted by $(F^k)_{k\in \N ^*}$. For proving that $f$ is a mild solution of \eqref{discrete-general-a}-\eqref{discrete-general-b}, it is sufficient to prove that for any $\eta >0$ and $i\in \{ 1, \cdot \cdot \cdot , p\} $, there is a set $X_i^\eta $ of $i$-characteristics with  complementary set of measure smaller than $c\eta $, such that 
\begin{align}\label{use-of-test-function}
\int _\Omega \varphi \chi _i^\eta f_i(z)dz&= \int _\Omega \varphi \chi _i^\eta f_{bi}(z_i^+(z))dz\nonumber \\
&+\int _\Omega \int _{-s_i^+(z)}^0\big( \varphi \chi _i^\eta Q_i(f,f)+\chi _i^\eta f_i\hspace*{0.02in}v_i\cdot \nabla \varphi \big) (z+sv_i)ds\hspace*{0.02in}dz,\quad \varphi \in C^1(\bar{\Omega }) ,
\end{align}
where $\chi _i^\eta $ denotes the characteristic function of $X_i^\eta $.
Define the set $X_i^\eta $ as follows. For every $\epsilon >0$, pass to the limit when $k\rightarrow +\infty $ in
\begin{align}\label{exp-form-chi-F-k}
&\chi _i^{k,\epsilon }F_i^k (z)\leq \chi _i^{k,\epsilon }F_i^k (z_i^-(z))\hspace*{0.02in}\exp \Big( \int _{-s_i^+(z)}^{s_i^-(z)}\nu _i^k (z+sv_i)ds\Big) , \quad \text{a.a.   }z\in \Omega , \quad k\in \N ^*,
\end{align}
and use the weak $L^1$ compactness of $(\chi _i^{k,\epsilon }F_i^k)_{k\in \N ^*}$, the weak $L^1$ compactness and the uniform boundedness in $L^\infty $ of $(\chi _i^{k,\epsilon }F_i^k(z_i^-(z)))_{k\in \N ^*}$, and the strong $L^1$ compactness of $(\int _{-s_i^+(z)}^{s_i^-(z)}\nu_i^k(z+sv_i)ds)_{k\in \N ^*}$. It implies that
\begin{align*}
&F_i^\epsilon (z)\leq F_i^\epsilon (z_i^-(z))\hspace*{0.02in}\exp \Big( \int _{-s_i^+(z)}^{s_i^-(z)}\nu _i(f)(z+sv_i)ds\Big) , \quad \text{a.a.   }z\in \Omega , \quad \epsilon \in ]0,1[,
\end{align*}
where $F_i^\epsilon $ is the limit of a subsequence of $(\chi _i^{k,\epsilon }F_i^k)_{k\in \N ^*}$ and $\nu _i(f)= \sum _{j,l,m=1}^p\Gamma _{ij}^{lm}f_j$. By the monotonicity in $\epsilon $ of $(F^\epsilon )_{\epsilon \in ]0,1[}$ (resp. $\big( F^\epsilon (z_i^-(z))\big) _{\epsilon \in ]0,1[})$ and the uniform boundedness of their masses, it holds that
\begin{align*}
&f_i(z)\leq f_i(z_i^-(z))\hspace*{0.02in}\exp \Big( \int _{-s_i^+(z)}^{s_i^-(z)}\nu _i(f)(z+sv_i)ds\Big) , \quad \text{a.a.   }z\in \Omega .
\end{align*}
From here the proof follows the lines of the proof of Lemma \ref{lemma-df-Omega-epsilon}, so that given $\eta >0$, there is a set $X_i^\eta $ of $i$-characteristics, with complementary set of measure smaller than $c\eta $, such that
\begin{align}\label{bdd-with-chi-eta}
&f_i(z)\leq \frac{1}{\eta }e^{\frac{p\Gamma }{\eta }}\quad \text{and}\quad \int _{-s_i^+(z)}^{s_i^-(z)}\nu _i(f)(z+sv_i)ds\leq \frac{p\Gamma }{\eta },\quad \text{a.a.  }z\in X_\eta .
\end{align}
Denote by $C^1_+(\bar{\Omega })$ the subspace of non-negative functions of $C^1(\bar{\Omega })$.
\hspace*{1.in}\\
%
%
% Lemma 4.1
%
\begin{lemma}\label{passage-limit-in-k}
\hspace*{1.in}\\
$f$ is a subsolution of \eqref{discrete-general-a}-\eqref{discrete-general-b}, i.e.
\begin{align}\label{f-subsolution}
\int _\Omega \varphi \chi _i^\eta f_i(z)dz&\leq  \int _\Omega \varphi f_{bi}(z_i^+(z))dz+\int _\Omega \int _{-s_i^+(z)}^0\chi _i^\eta f_i\hspace*{0.02in}v_i\cdot \nabla \varphi  (z+sv_i)ds\hspace*{0.02in}dz\nonumber \\
&+\int _\Omega \int _{-s_i^+(z)}^0\varphi \hspace*{0.01in}Q_i(f,f) (z+sv_i)ds\hspace*{0.02in}dz,\quad 1\leq i\leq p,\quad \varphi \in C^1_+(\bar{\Omega }).
\end{align}
\end{lemma}
\underline{Proof of Lemma \ref{passage-limit-in-k}.}\\
Let $i\in \{ 1,\cdot \cdot \cdot ,p\} $ and $\varphi \in C^1_+(\bar{\Omega })$ be given. Write the mild form of $\varphi \chi _i^\eta \chi _i^{k,\epsilon }F_i^k$ and integrate it on $\Omega $. It results
\begin{align}\label{passage-limit-1}
\int _\Omega \varphi \chi _i^\eta \chi _i^{k,\epsilon }F_i^k(z)dz&= \int _\Omega \varphi \chi _i^\eta \chi _i^{k,\epsilon }f_{bi}^k(z_i^+(z))dz+\int _\Omega \int _{-s_i^+(z)}^0\chi _i^\eta \chi _i^{k,\epsilon }F_i^k\hspace*{0.02in}v_i\cdot \nabla \varphi (z+sv_i)ds\hspace*{0.02in}dz\nonumber \\
&+\int _\Omega \int _{-s_i^+(z)}^0\varphi \chi _i^\eta \chi _i^{k,\epsilon }\big( Q_i^{+k}-F_i^k\nu _i^k\big) (z+sv_i)ds\hspace*{0.02in}dz.
\end{align}
By the weak $L^1$ compactness of $(F_i^k)_{k\in \N ^*}$ and the linearity with respect to $\chi _i^{k,\epsilon }F_i^k$ of the first line of \eqref{passage-limit-1}, its passage to the limit when $k\rightarrow +\infty $ is straightforward. Let us pass to the limit when $k\rightarrow +\infty $ in any term of the loss term of \eqref{passage-limit-1}, denoted by $\Gamma _{ij}^{lm}L^k$, where 
\begin{align}\label{passage-limit-2}
&L^k:= \int _\Omega \chi _i^\eta \chi _i^{k,\epsilon }(z)\int _{-s_i^+(z)}^0\varphi \frac{F_i^k}{1+\frac{F_i^k}{k}}\frac{F_j^k}{1+\frac{F_j^k}{k}}(z+sv_i)ds\hspace*{0.02in}dz,\quad j\in J_i,
\end{align}
and $J_i$ is defined in \eqref{df-J-i}.
By integration by parts, $L_k$ equals
\begin{align}\label{passage-limit-3}
&\int _\Omega \int _{-s_i^+(z)}^0\chi _i^\eta \chi _i^{k,\epsilon }\big( \varphi (Q_i^{+k}-F_i^k\nu _i^k) +(v_i\cdot \nabla \varphi )F_i^k\big)(z+sv_i)\Big( \int _{s}^0\chi _i^{k,\epsilon }\frac{F_j^k}{(1+\frac{F_i^k}{k})(1+\frac{F_j^k}{k})}(z+rv_i)dr\Big) ds\hspace*{0.02in}dz\nonumber \\
&+\int _\Omega \chi _i^\eta \chi _i^{k,\epsilon }\varphi \frac{f_{bi}^k}{1+\frac{f_{bi}^k}{k}}(z_i^+(z))\int _{-s_i^+(z)}^0\frac{F_j^k}{1+\frac{F_j^k}{k}}(z+sv_i)ds\hspace*{0.02in}dz.
\end{align}
Denote by $(a,b)$ the coordinate system in the $(v_i,v_j)$ basis, $(a^-, a^+)\in \R ^2$ and $(b^-(a), b^+(a))\in \R ^2$ for every $a\in ]a^-, a^+[ $, such that
\begin{align}\label{decomposition-Omega}
&\Omega = \{ av_i+bv_j;\quad a\in ]a^-, a^+[, \quad b\in ]b^-(a), b^+(a)[ \hspace*{0.02in}\} .
\end{align}
The first term in $L^k$ can be written as $\int _{a^-}^{a^+}l^k(a)da$ with $l^k$ defined as
\begin{align}\label{passage-limit-5}
&l^k(a)= \int _{b^- (a)}^{b^+(a)}\int _{-s_i(bv_j)}^a\chi _i^\eta \chi _i^{k,\epsilon }\big( \varphi (Q_i^{+k}-F_i^k\nu _i^k)+(v_i\cdot \nabla \varphi )F_i^k\big) (sv_i+bv_j)\nonumber \\
&\hspace*{1.8in}\big( \int _{s}^a\chi _i^{k,\epsilon }\frac{F_j^k}{(1+\frac{F_i^k}{k})(1+\frac{F_j^k}{k})}(rv_i+bv_j)dr\big) ds\hspace*{0.02in}db.
\end{align}
For each rational number $a$, the sequence of functions 
\begin{align*}
&(b,s)\in [b^-(a),b^+(a)]\times [-s_i^+(bv_j),a]\rightarrow \chi _i^\eta \chi _i^{k,\epsilon }\big( \varphi (Q_i^{+k}-F_i^k\nu _i^k)+(v_i\cdot \nabla \varphi )F_i^k\big) (sv_i+bv_j)
\end{align*}
 is weakly compact in $L^1$, whereas 
\begin{align*} 
&(b,s)\rightarrow \int _{s}^a\chi _i^{k,\epsilon }\frac{F_j^k}{(1+\frac{F_i^k}{k})(1+\frac{F_j^k}{k})}(rv_i+bv_j)dr
\end{align*}
is by Lemma \ref{compactness-integrated-collision-frequency} strongly compact in $L^1$, and by Lemma \ref{lemma-df-Omega-epsilon} uniformly bounded in $L^\infty $. The convergence follows for any rational number $a$. With a diagonal process, there is a subsequence of $(l^k)$, still denoted by $(l^k)$, converging for any rational $a$. Moreover,
\begin{align}\label{passage-limit-6}
&\lim _{h\rightarrow 0}\hspace*{0.02in}\big( l^k(a+h)-l^k(a)\big) = 0,
\end{align}
uniformly with respect to $k$ and $a$, by the weak $L^1$ compactness of 
\begin{align*}
&\big( \chi _i^\eta \chi _i^{k,\epsilon }( \varphi (Q_i^{+k}-F_i^k\nu _i^k)+(v_i\cdot \nabla \varphi )F_i^k\big) _{k\in \N ^*}\quad \text{and}\quad (F_j^k)_{k\in \N ^*}.
\end{align*}
Thus $(l^k)$ is a uniform converging sequence on $[a^-, a^+]$. The second term in $L^k$ can be treated analogously, $(\chi _i^{k,\epsilon }f_{bi}^k)_{k\in \N ^*}$ being uniformly bounded in $L^\infty $. The convergence follows.\\
\hspace*{1.in}\\
In order to determine the limit of $L^k$ when $k\rightarrow +\infty $, remark that 
\begin{align*}
&\chi _i^\eta \chi _i^{k,\epsilon }( \varphi (Q_i^{+k}-F_i^k\nu _i^k)+(v_i\cdot \nabla \varphi )F_i^k= v_i\cdot \nabla (\chi _i^\eta \chi _i^{k,\epsilon }\varphi \hspace*{0.01in}F_i^k),
\end{align*} 
which weakly converges in $L^1$ to $v_i\cdot \nabla (\chi _i^\eta \varphi \hspace*{0.01in}F_i^\epsilon )$ when $k\rightarrow +\infty $. Hence
\begin{align}
\lim _{k\rightarrow +\infty }L^k&=\int _\Omega \int _{-s_i^+(z)}^0v_i\cdot \nabla (\chi _i^\eta \varphi \hspace*{0.01in}F_i^\epsilon )(z+sv_i)\Big( \int _{s}^0f_j(z+rv_i)dr\Big) ds\hspace*{0.02in}dz\nonumber \\
&+\int _\Omega \chi _i^\eta \varphi f_{bi}(z_i^+(z))\Big( \int _{-s_i^+(z)}^0f_j(z+sv_i)ds\Big) \hspace*{0.02in}dz.\nonumber 
\end{align}
By a backwards integration by parts,
\begin{align}\label{passage-limit-8}
\lim _{k\rightarrow +\infty }L^k&= \int _\Omega \int _{-s_i^+(z)}^0\varphi \hspace*{0.01in}\chi _i^\eta F_i^\epsilon f_j(z+sv_i)ds\hspace*{0.02in}dz.
\end{align} 
\hspace*{0.02in}\\
In order to prove \eqref{f-subsolution}, let us prove that each 
\begin{align}\label{passage-limit-9}
&\Gamma _{ij}^{lm}\int _\Omega \int _{-s_i^+(z)}^0\varphi \chi _i^\eta \chi _i ^{k,\epsilon }\frac{F_l^k}{1+\frac{F_l^k}{k}}\frac{F_m^k}{1+\frac{F_m^k}{k}}(z+sv_i)ds\hspace*{0.02in}dz,\quad j\in J_i,
\end{align}
term from $Q_i^{+k}$ in \eqref{passage-limit-1} converges when $k\rightarrow +\infty $ to a limit smaller than 
\begin{align}\label{passage-limit-10}
&\Gamma _{ij}^{lm}\int _\Omega \int _{-s_i^+(z)}^0\varphi \chi _i^\eta F_l^{\epsilon ^\prime }f_m(z+sv_i)ds\hspace*{0.02in}dz+\alpha (\epsilon ^\prime ),\hspace*{0.04in}\epsilon ^\prime \in ] 0,1[ ,\quad \text{with  } \lim _{\epsilon ^\prime \rightarrow 0}\alpha (\epsilon ^\prime )= 0.
\end{align} 
Take $\Gamma _{ij}^{lm}= 1$, $j\in J_i$, for simplicity.  $(\mu _ {\frac{1}{n}})_{n\in \N ^*}$ being the sequence of mollifiers defined at the beginning of Section \ref{Section2} for $\alpha = \frac{1}{n}$, split \eqref{passage-limit-9} into
\begin{align}\label{passage-limit-12}
&\int _\Omega \int _{-s_i^+(z)}^0\varphi (\chi _i^\eta \ast \mu _{\frac{1}{n}})\chi _l^{k,\epsilon ^\prime }\chi _i^{k,\epsilon }\frac{F_l^k}{1+\frac{F_l^k}{k}}\frac{F_m^k}{1+\frac{F_m^k}{k}}(z+sv_i)ds\hspace*{0.02in}dz\nonumber \\
&+\int _\Omega \int _{-s_i^+(z)}^0\varphi (\chi _i^\eta \ast \mu _{\frac{1}{n}})(1-\chi _l^{k,\epsilon ^\prime })\chi _i^{k,\epsilon }\frac{F_l^k}{1+\frac{F_l^k}{k}}\frac{F_m^k}{1+\frac{F_m^k}{k}}(z+sv_i)ds\hspace*{0.02in}dz\nonumber \\
&\int _\Omega \int _{-s_i^+(z)}^0\varphi \big( \chi _i^\eta -(\chi _i^\eta \ast \mu _{\frac{1}{n}})\big) \chi _i^{k,\epsilon }\frac{F_l^k}{1+\frac{F_l^k}{k}}\frac{F_m^k}{1+\frac{F_m^k}{k}}(z+sv_i)ds\hspace*{0.02in}dz\nonumber \\
&\hspace*{1.in}\nonumber \\
&\leq \int _\Omega \int _{-s_i^+(z)}^0\varphi (\chi _i^\eta \ast \mu _{\frac{1}{n}})\chi _l^{k,\epsilon ^\prime }\frac{F_l^k}{1+\frac{F_l^k}{k}}\frac{F_m^k}{1+\frac{F_m^k}{k}}(z+sv_i)ds\hspace*{0.02in}dz\nonumber \\
&+\frac{c}{\ln \Lambda }+\frac{c\Lambda }{\epsilon ^2}e^{\frac{p\Gamma }{\epsilon ^2}}\sum _{j\in J_i}\Big( \int _{\Omega _l^{k,\epsilon ^\prime }}F_j^k(z)dz+\int _\Omega \varphi \mid \chi _i^\eta -(\chi _i^\eta \ast \mu _{\frac{1}{n}})\mid F_j^k(z)dz \Big) \nonumber \\
&\hspace*{1.in}\nonumber \\
&\leq \int _\Omega \int _{-s_i^+(z)}^0\varphi (\chi _i^\eta \ast \mu _{\frac{1}{n}})\chi _l^{k,\epsilon ^\prime }\frac{F_l^k}{1+\frac{F_l^k}{k}}\frac{F_m^k}{1+\frac{F_m^k}{k}}(z+sv_i)ds\hspace*{0.02in}dz\nonumber \\
&\hspace*{1.in}\nonumber \\
&+\frac{c}{\ln \Lambda }+\frac{c\Lambda }{\epsilon ^2}e^{\frac{p\Gamma }{\epsilon ^2}}\Big( \Lambda ^\prime \epsilon ^\prime  +\frac{1}{\ln \Lambda ^\prime }+\frac{1}{\ln \frac{k}{2}}+\tilde{\Lambda }\parallel \chi _i^\eta -(\chi _i^\eta \ast \mu _{\frac{1}{n}})\parallel _{L^1}+\frac{1}{\ln \tilde{\Lambda }}\Big) ,\quad \text{by   }\eqref{weak-compactness-1},\nonumber \\
&\hspace*{3.in}\Lambda >1,\quad \Lambda ^\prime >1,\quad \tilde{\Lambda }>1,\quad \epsilon ^\prime >0.
\end{align}
Denote by $D$ the Jacobian of the change of variables $z\rightarrow (a,b)$. For some smooth function $A$, and any integrable function $g$, 
\begin{align*}
\int _\Omega \int _{-s_i^+(z)}^0g(z+sv_i)dsdz&=  D\int _{b^-}^{b^+}\int _{a^-(b)}^{a^+(b)}\int _{-s_i^+(bv_j)}^ag(sv_i+bv_j)ds\hspace*{0.02in}da\hspace*{0.02in}db\\
&=  D\int _{b^-}^{b^+}\int _{-s_i^+(bv_j)}^{a^+(b)}(a^+(b)-\max \{ a^-(b), s\} )g(sv_i+bv_j)\hspace*{0.02in}ds\hspace*{0.02in}db\\
&=  \int _{\Omega }A(\alpha ,\gamma )g(\alpha v_l+\gamma v_m)\hspace*{0.02in}d\alpha \hspace*{0.02in}d\gamma .
\end{align*}
Hence, 
\begin{align}
&\lim _{k\rightarrow +\infty }\int \int _{-s_i^+(z)}^0\varphi (\chi _i^\eta \ast \mu _{\frac{1}{n}})\chi _l^{k,\epsilon ^\prime }\frac{F_l^k}{1+\frac{F_l^k}{k}}\frac{F_m^k}{1+\frac{F_m^k}{k}}(z+sv_i)ds\hspace*{0.02in}dz\nonumber \\
&= \int _\Omega \int _{-s_i^+(z)}^0\varphi (\chi _i^\eta \ast \mu _{\frac{1}{n}})F_l^{\epsilon ^\prime }f_m(z+sv_i)ds\hspace*{0.02in}dz,\quad \epsilon ^\prime \in ] 0,1[ .
\end{align}
For $\tilde{\Lambda }$ large enough, pass to the limit when $k\rightarrow +\infty $ and $n\rightarrow +\infty $ in \eqref{passage-limit-12}.  Up to subsequences, the weak $L^1$ limits $F_i^\epsilon $ and $F_i^{\epsilon ^\prime }$ of $(\chi _i^{k,\epsilon }F_i^k)_{k\in \N ^*}$ and $(\chi _i^{k,\epsilon ^\prime }F_i^k)_{k\in \N ^*}$ when $k\rightarrow +\infty $ satisfy 
\begin{align}\label{passage-limit-a}
&\int _\Omega \varphi \chi _i^\eta F_i^\epsilon (z)dz\leq \int _\Omega \varphi \chi _i^\eta f_{bi}^k(z_i^+(z))dz+\int _\Omega \int _{-s_i^+(z)}^0\chi _i^\eta F_i^\epsilon \hspace*{0.02in}v_i\cdot \nabla \varphi (z+sv_i)ds\hspace*{0.02in}dz\nonumber \\
&+\int _\Omega \int _{-s_i^+(z)}^0\varphi \chi _i^\eta \big( Q_i^+(F^{\epsilon ^\prime },f)-F_i^\epsilon \nu _i(f)\big) (z+sv_i)ds\hspace*{0.02in}dz\\
&+\frac{c}{\ln \Lambda }+\frac{c\Lambda }{\epsilon ^2}e^{\frac{p\Gamma }{\epsilon ^2}}\big( \Lambda ^\prime \epsilon ^\prime  +\frac{1}{\ln \Lambda ^\prime }\big) ,\quad (\epsilon ,\epsilon ^\prime )\in ]0,1[ ^2 ,\quad \Lambda >1, \quad \Lambda ^\prime >1.\nonumber 
\end{align}
Choose $\Lambda $ large enough, $\epsilon $ small enough, $\Lambda ^\prime $ large enough,  $\epsilon ^\prime $ small enough, in this order.
The passage to the limit when $\epsilon \rightarrow 0$ and $\epsilon ^\prime \rightarrow 0$ in \eqref{passage-limit-a} results from the monotone convergence theorem, the family $(F^\epsilon )_{\epsilon \in ]0,1[ }$ being non decreasing, with mass uniformly bounded, together with the mass of $(\chi _i^\eta Q_i^+(F^{\epsilon ^\prime },f))_{\epsilon ^\prime \in ]0,1[ }$ and $(\chi _i^\eta F_i^{\epsilon ^\prime }\nu _i(f))_{\epsilon ^\prime \in ]0,1[ }$. Consequently, \eqref{f-subsolution} holds. \cqfd
%
%
%  Lemma 4.2
%
%
\begin{lemma}\label{proof-theorem-1-1}
$f$ is a solution to \eqref{discrete-general-a}-\eqref{discrete-general-b}.
\end{lemma}
\underline{Proof of Lemma \ref{proof-theorem-1-1}.}\\
For proving Lemma \ref{proof-theorem-1-1}, it remains to prove that 
\begin{align}\label{passage-limit-18}
\int _\Omega \varphi \chi _i^\eta f_i(z)dz&\geq  \int _\Omega \varphi \chi _i^\eta f_{bi}(z_i^+(z))dz+\int _\Omega \int _{-s_i^+(z)}^0\chi _i^\eta f_i\hspace*{0.02in}v_i\cdot \nabla \varphi (z+sv_i)ds\hspace*{0.02in}dz\nonumber \\
&+\int _\Omega \int _{-s_i^+(z)}^0\varphi \hspace*{0.02in}\chi _i^\eta Q_i(f,f)(z+sv_i)ds\hspace*{0.02in}dz,\quad 1\leq i\leq p,\quad \varphi \in C^1_+(\bar{\Omega }).
\end{align}
For $\beta>0$, start from the equation for $\varphi \chi _i^\eta F_i^k$ written in renormalized form, 
\begin{align}\label{passage-limit-19}
&\beta^{-1}\varphi \chi _i^\eta \ln (1+\beta F_i^k) (z)-\beta^{-1} \varphi \chi _i^\eta \ln (1+\beta f_{bi}^k)(z_i^+(z))\nonumber \\
&+\int _{-s_i^+(z)}^0\beta^{-1}\chi _i^\eta \ln (1+\beta F_i^k) \hspace*{0.02in}v_i\cdot \nabla \varphi (z+sv_i)ds= \int _{-s_i^+(z)}^0\frac{\varphi  \hspace*{0.02in}\chi _i^\eta (Q_i^{+k}-F_i^k\nu _i^k)}{1+\beta F_i^k}(z+sv_i)ds.
\end{align} 
It holds 
\begin{align*}
&\beta ^{-1}\ln(1+\beta x)<x,\hspace*{0.04in} \beta \in ]0,1[ \quad \text{  and  }\quad \lim_{\beta \rightarrow 0} \beta ^{-1}\ln(1+\beta  x)= x,\quad x>0.
\end{align*} 
Hence in weak $L^1$  the sequence $(\beta ^{-1}\ln \big( 1+\beta F^k_i\big) )_{k\in \N ^*}$ converges modulo subsequence to a function $F^\beta \leq f$ 
  when $k\rightarrow +\infty $. The mass of the limit increases to the mass  of $f$, when $\beta
 \rightarrow 0$. 
 This gives in the final limit $\beta \rightarrow 0$ for the l.h.s. of \eqref{passage-limit-19},
 \begin{align}\label{passage-limit-19-a}
&\varphi \chi _i^\eta f_i(z)- \varphi \chi _i^\eta f_{bi} (z_i^+(z))
- \int _{-s_i^+(z)}^0 \chi _i^\eta f_i \hspace*{0.02in}v_i\cdot \nabla \varphi (z+sv_i)ds.
\end{align}
Using analogous arguments as for the limit of the loss term in Lemma \ref{passage-limit-in-k}, it holds that
\begin{align*}
&\lim _{k\rightarrow +\infty }\Gamma _{ij}^{lm}\int _\Omega \int _{-s_i^+(z)}^0\frac{\varphi \chi _i^\eta F_i^kF_j^k}{1+\beta F_i^k}(z+sv_i)ds\hspace*{0.02in}dz\nonumber \\
&= \Gamma _{ij}^{lm}\int _\Omega \int _{-s_i^+(z)}^0\varphi \chi _i^\eta \big( \underset{k\rightarrow +\infty }{\text{weak}L^1\text{lim}}\frac{F_i^k}{1+\beta F_i^k}\big) f_j(z+sv_i)ds\hspace*{0.02in}dz,\quad j\in J_i.
\end{align*}
But 
\begin{align*}
&\underset{k\rightarrow +\infty }{\text{weak}L^1\text{lim}}\frac{F_i^k}{1+\beta F_i^k}\leq \underset{k\rightarrow +\infty }{\text{weak}L^1\text{lim}}F_i^k,
\end{align*}
and
\begin{align*}
&\int _\Omega \underset{k\rightarrow +\infty }{\text{weak}L^1\text{lim}}\frac{F_i^k}{1+\beta F_i^k}(z)dz \quad \text{increases to     }\int _\Omega \underset{k\rightarrow +\infty }{\text{weak}L^1\text{lim}}F_i^k(z)dz
\end{align*}
when $\beta \rightarrow 0$. Hence
\begin{align}\label{passage-limit-}
&\lim _{\beta \rightarrow 0}\lim _{k\rightarrow +\infty }\Gamma _{ij}^{lm}\int _\Omega \int _{-s_i^+(z)}^0\frac{\varphi \chi _i^\eta F_i^kF_j^k}{1+\beta F_i^k}(z+sv_i)ds\hspace*{0.02in}dz= \Gamma _{ij}^{lm}\int _\Omega \int _{-s_i^+(z)}^0\varphi \chi _i^\eta f_if_j(z+sv_i)ds\hspace*{0.02in}dz.
\end{align}
For the gain term and any $(l,m)\in \{ 1, \cdot \cdot \cdot , p\} ^2$ such that $\Gamma _{ij}^{lm}>0$ for some $j\in \{ 1,\cdot \cdot \cdot , p\} $,
\begin{align}\label{limit-in-k-24}
&\int _\Omega \int _{-s_i^+(z)}^0\frac{\varphi \chi _i^\eta }{1+\beta F_i^k}\frac{F_l^k}{1+\frac{F_l^k}{k}}\frac{F_m^k}{1+\frac{F_m^k}{k}}(z+sv_i)ds\hspace*{0.02in}dz\nonumber \\
&\geq \int _\Omega \int _{-s_i^+(z)}^0\frac{\varphi \chi _i^\eta \chi _l^{k,\epsilon }}{1+\beta F_i^k}\frac{F_l^k}{1+\frac{F_l^k}{k}}\frac{F_m^k}{1+\frac{F_m^k}{k}}(z+sv_i)ds\hspace*{0.02in}dz\nonumber \\
&= \int _\Omega \int _{-s_i^+(z)}^0\varphi \chi _i^\eta \chi _l^{k,\epsilon }\frac{F_l^k}{1+\frac{F_l^k}{k}}\frac{F_m^k}{1+\frac{F_m^k}{k}}(z+sv_i)ds\hspace*{0.02in}dz\nonumber \\
&-\int _\Omega \int _{-s_i^+(z)}^0\varphi \chi _i^\eta \chi _l^{k,\epsilon }\frac{\beta F_i^k}{1+\beta F_i^k}\frac{F_l^k}{1+\frac{F_l^k}{k}}\frac{F_m^k}{1+\frac{F_m^k}{k}}(z+sv_i)ds\hspace*{0.02in}dz\nonumber \\
&\geq \int _\Omega \int _{-s_i^+(z)}^0\varphi \chi _i^\eta \chi _l^{k,\epsilon }\frac{F_l^k}{1+\frac{F_l^k}{k}}\frac{F_m^k}{1+\frac{F_m^k}{k}}(z+sv_i)ds\hspace*{0.02in}dz\nonumber \\
&-c\Lambda \sum _{j\in J_i}\int _\Omega \int _{-s_i^+(z)}^0\varphi \chi _i^\eta \chi _l^{k,\epsilon }\frac{\beta (F_i^k)^2F_j^k}{1+\beta F_i^k}(z+sv_i)ds\hspace*{0.02in}dz-\frac{c}{\ln \Lambda }\quad \Lambda >1,\quad \epsilon \in ]0,1[.
\end{align}
It holds
\begin{align}\label{limit-in-k-25}
&\lim _{k\rightarrow +\infty }\int _\Omega \int _{-s_i^+(z)}^0\varphi \chi _i^\eta \chi _l^{k,\epsilon }\frac{F_l^k}{1+\frac{F_l^k}{k}}\frac{F_m^k}{1+\frac{F_m^k}{k}}(z+sv_i)ds\hspace*{0.02in}dz= \int _\Omega \int _{-s_i^+(z)}^0\varphi \chi _i^\eta F_l^\epsilon  f_m^k(z+sv_i)ds\hspace*{0.02in}dz.
\end{align}
Choose $\Lambda $ large enough and split the domain of integration of every $j\in J_i$ term in \eqref{limit-in-k-24} into 
\begin{align*}
\{ F_i^k\leq \Lambda ^\prime  \} &\cup  \{ F_i^k>\Lambda ^\prime  \quad \text{and}\quad F_i^kF_j^k> \tilde{\Lambda } \frac{F_l^k}{1+\frac{F_l^k}{k}}\frac{F_m^k}{1+\frac{F_m^k}{k}}\} \\
&\cup \{ F_i^k>\Lambda ^\prime \quad \text{and}\quad F_i^kF_j^k\leq \tilde{\Lambda } \frac{F_l^k}{1+\frac{F_l^k}{k}}\frac{F_m^k}{1+\frac{F_m^k}{k}}\} ,\quad \Lambda ^\prime  >1,\quad \tilde{\Lambda } >1.
\end{align*}
It holds that
\begin{align}\label{limit-in-k-24-a}
&\int _\Omega \int _{-s_i^+(z)}^0\varphi \chi _i^\eta \chi _l^{k,\epsilon }\frac{\beta (F_i^k)^2F_j^k}{1+\beta F_i^k}(z+sv_i)ds\hspace*{0.01in}dz\leq c\Big( \beta (\Lambda ^\prime )^2+\frac{1}{\ln \tilde{\Lambda }}+\frac{\tilde{\Lambda }}{\epsilon ^2}e^{\frac{p\Gamma }{\epsilon ^2}}\int _{F_i^k>\Lambda ^\prime }F_m^k(z)dz\Big) ,\nonumber \\
&\hspace*{3.in} \beta \in ]0,1[ ,\quad \Lambda ^\prime >0,\quad \tilde{\Lambda }>1.
\end{align}
The last term in \eqref{limit-in-k-24-a} tends to zero when $\tilde{\Lambda }\rightarrow +\infty $, $\Lambda ^\prime \rightarrow +\infty $, $\beta \rightarrow 0$ in this order, uniformly with respect to $k$. Consequently,
\begin{align*}
&\lim _{\beta \rightarrow 0}\lim _{k\rightarrow +\infty }\int _\Omega \int _{-s_i^+(z)}^0\frac{\varphi \chi _i^\eta }{1+\beta F_i^k}\frac{F_l^k}{1+\frac{F_l^k}{k}}\frac{F_m^k}{1+\frac{F_m^k}{k}}(z+sv_i)ds\hspace*{0.01in}dz\geq \int _\Omega \int _{-s_i^+(z)}^0\varphi \chi _i^\eta F_l^\epsilon f_m(z+sv_i)ds\hspace*{0.01in}dz.
\end{align*}
This holds for every $\epsilon >0$. Hence
\begin{align}\label{limit-in-k-24-d}
&\lim _{\beta \rightarrow 0}\lim _{k\rightarrow +\infty }\int _\Omega \int _{-s_i^+(z)}^0\frac{\varphi \chi _i^\eta }{1+\beta F_i^k}\frac{F_l^k}{1+\frac{F_l^k}{k}}\frac{F_m^k}{1+\frac{F_m^k}{k}}(z+sv_i)ds\hspace*{0.01in}dz\geq \int _\Omega \int _{-s_i^+(z)}^0\varphi \chi _i^\eta f_l f_m(z+sv_i)ds\hspace*{0.01in}dz.
\end{align}
And so, \eqref{passage-limit-18} holds. Together with \eqref{f-subsolution}, this proves \eqref{use-of-test-function}. \cqfd

\end{document}